\documentclass[11pt,a4paper]{amsart} 

\usepackage{latexsym}
\usepackage{geometry}
\usepackage{pinlabel}         
\usepackage{graphicx}
\usepackage{amssymb, amsthm}
\usepackage{epstopdf}
\usepackage{romannum}
\usepackage{tikz}
\usepackage{subfig}
\usepackage{MnSymbol}
\usepackage[arrow, matrix, curve]{xy} 

\theoremstyle{plain}
\newtheorem*{NewTheorem}{Theorem 1}
\newtheorem*{definition}{Definition}
\newtheorem*{remark}{Remark}

\title[]{On property (T) for automorphism groups of graph products}
\author{Nils Leder and Olga Varghese}
\thanks{Funded by the Deutsche 
Forschungsgemeinschaft (DFG, German Research Foundation) under Germany's 
Excellence Strategy -EXC 2044-, Mathematics M\"unster: Dynamics-Geometry-Structure}
\date{\today}

\address{Nils Leder\\
Department of Mathematics\\
M\"unster University\\ 
Einsteinstra\ss e 62\\
48149 M\"unster (Germany)}
\email{n\_lede02@uni-muenster.de}

\address{Olga Varghese\\
Department of Mathematics\\
M\"unster University\\ 
Einsteinstra\ss e 62\\
48149 M\"unster (Germany)}
\email{olga.varghese@uni-muenster.de}

\begin{document}
\pagenumbering{arabic}
\begin{abstract}
We show that the automorphism group of a graph product of finite groups ${\rm Aut}(G_\Gamma)$ has Kazhdan's property (T) if and only if $\Gamma$ is a complete graph.
\end{abstract}
\maketitle

\section*{Property (T) for ${\rm Aut}(G_\Gamma)$}
Property (T) was defined by Kazhdan for locally compact groups in terms of unitary representations. This property was reformulated in many different mathematical areas, in particular in geometric group theory. It was proven by Delorme and Guichardet that a countable group $G$ has property (T) if and only if every action of $G$ on a Hilbert space by isometries has a global fixed point, see (\cite[2.12.4]{Bekka}). Examples of non-finite groups with property (T) are the general linear groups ${\rm GL}_n(\mathbb{Z})$ for $n\geq 3$ (\cite[4.2.5]{Bekka}). The abelianization  of a free group $F_n$ induces an epimorphism ${\rm Aut}(F_n)\twoheadrightarrow{\rm GL}_n(\mathbb{Z})$. Thus, it is natural to ask if the group ${\rm Aut}(F_n)$ for $n\geq 3$ also has property (T). It was proven by McCool in \cite{McCool} that ${\rm Aut}(F_3)$ does not have property (T). Further, Kaluba, Nowak and Ozawa proved in \cite{Kaluba} that ${\rm Aut}(F_5)$ has property (T) and recently, Kaluba, Kielak and Nowak  showed in \cite{Kielak} that 
${\rm Aut}(F_n)$ for $n\geq 6$ has property (T). It is still an open question whether ${\rm Aut}(F_4)$ has property (T) or not. 

The free group $F_n$ is an example of a graph product where the graph has $n$ vertices, no edges and each vertex is labeled with the infinite cyclic group $\mathbb{Z}$. More precisely, let ${\Gamma=(V, E)}$ be a finite simplicial graph with vertex labeling ${V\rightarrow\left\{\text{non-trivial groups}\right\}, v\mapsto G_v}$. The graph product $G_\Gamma$ is the group obtained from the free product of the $G_v$, $v\in V$, by adding the relations $gh=hg$ for all $g\in G_v, h\in G_w$ such that $\left\{v, w\right\}\in E$. 
We are interested in graph products of finite groups, i. e. for all $v\in V$ the vertex group $G_v$ is finite. Our leading question is the following:
\begin{center}
{\it For which shape of the graph $\Gamma$ does the automorphism group ${\rm Aut}(G_\Gamma)$ have property (T)?}
\end{center} There are some partial results in the literature. It was proven in \cite{Varghese} that the group ${\rm Aut}(\mathbb{Z}/2*\ldots*\mathbb{Z}/2)$ does not have property (T). This result was generalized in \cite{Leder} for the groups ${\rm Aut}(\mathbb{Z}/n_1*\ldots*\mathbb{Z}/n_m)$. Further, there are results regarding outer automorpism groups of graph products of finite abelian groups in \cite[Cor. 3]{Sale}.
\newpage
We prove the following characterization:
\begin{NewTheorem}
Let $G_\Gamma$ be a graph product of finite groups. Then  ${\rm Aut}(G_\Gamma)$ has property $(T)$ if and only if $\Gamma$ is a complete graph.
\end{NewTheorem}

For the proof we need to recall some important properties of graph products. The following definition is crucial in our argument.

\begin{definition}
For any $V' \subseteq V$ we denote by $G_{V'}$ the subgroup generated by all vertex groups $G_v, v \in V'$. We call $G_{V'}$ the \emph{special subgroup} defined by $V'$.
\end{definition}

Note that any special subgroup $G_{V'}$ is isomorphic to the graph product defined by the subgraph of $\Gamma$ spanned by $V'$.

\begin{remark}
The graph product $G_\Gamma$ has finitely many conjugacy classes of finite subgroups. This follows from the fact that each finite subgroup of $G_\Gamma$ lies in a parabolic subgroup $gG_{\Gamma'}g^{-1}$ where $\Gamma'$ is a complete subgraph of $\Gamma$, see \cite[Lemma 4.5]{Green}.
\end{remark}

The proof of Theorem 1 also involves Serre's fixed point property ${\rm F}\mathcal{A}$ which is closely connected to property (T).

\begin{definition}
A group $G$ is said to have \emph{property} ${\rm F}\mathcal{A}$ if for any action of $G$ on a simplicial tree $T$, without inversions of edges, there exists a global fixed point, i.e. there exists a vertex $v \in V(T)$ such that $g(v)=v$ for all $g \in G$.
\end{definition}

It was proven independently by Alperin in \cite{Alperin} and by Watatani in \cite{Watatani} that if a countable group $G$ has property (T), then
$G$ has Serre's property ${\rm F}\mathcal{A}$.

\begin{proof}[Proof of Theorem 1]
If $\Gamma$ is a complete graph, then the group $G_\Gamma$ is a direct product of finite groups $G_v$, $v\in V$. In particular, $G_\Gamma$ is finite. Thus ${\rm Aut}(G_\Gamma)$ is also finite and therefore has property (T) by Bruhat-Tits fixed point theorem \cite[II 2.8]{Bridson}. 

For the other implication of Theorem 1, assume that $\Gamma$ is not complete. Since $G_\Gamma$ is finitely generated, the group ${\rm Aut}(G_\Gamma)$ is countable. As property (T) descends to finite index subgroups (\cite[2.5.7]{Bekka}) it follows by Alperin and Watatani's result, that if a countable group $G$ has property (T), then any finite index subgroup of $G$ has property ${\rm F}\mathcal{A}$. Therefore it is sufficient to construct a finite index subgroup of ${\rm Aut}(G_\Gamma)$ that does not have property ${\rm F}\mathcal{A}$.
 
By the above remark, the group $G_\Gamma$ contains finitely many conjugacy classes of finite subgroups. Let $n$ be the cardinality of the set of conjugacy classes of finite subgroups of $G_\Gamma$. Then the group ${\rm Aut}(G_\Gamma)$ acts on this set in the canonical way. Thus, we obtain a homomorphism
\[
\Phi:{\rm Aut}(G_\Gamma)\rightarrow{\rm Sym}(n).
\]
By construction the kernel ${\rm ker}(\Phi)$ has finite index in ${\rm Aut}(G_\Gamma)$.
We claim that ${\rm ker}(\Phi)$ does not have Serre's property ${\rm F}\mathcal{A}$.

Since $\Gamma$ is not complete, there exist two vertices $v, w\in V$ of $\Gamma$ which are not connected by an edge. 
Let $\pi: G_\Gamma\twoheadrightarrow G_{\{v,w\}} \cong G_v*G_w$ be the projection defined as follows: $\pi(g)=g$ for $g\in G_v$, $\pi(g)=g$ for $g\in G_w$ and $\pi(g)=1$ for $g\in G_x$ for $x\in V-\left\{v,w\right\}$. The kernel of $\pi$ is the normal closure of the subgroup $G_{V-\left\{v,w\right\}}$ and is  characteristic under ${\rm ker}(\Phi)$. Hence we get the following canonical homomorphism 
\[
\Psi:{\rm ker}(\Phi)\rightarrow{\rm Aut}(G_v*G_w).
\]
The subgroup ${\rm Inn}(G_v*G_w)$ is contained in $\Psi({\rm ker}(\Phi))$. Further, it was proven in \cite{Karrass} and in \cite[Theorem 1]{Pettet} that the action of  $G_v*G_w\cong{\rm Inn}(G_v*G_w)$ on the associated Bass-Serre tree \cite[\S 4, Thm. 7]{Serre} extends to the whole automorphism group ${\rm Aut}(G_v*G_w)$, in particular to the subgroup $\Psi({\rm ker}(\Phi))$. Hence ${\rm ker}(\Phi)$ does not have property ${\rm F}\mathcal{A}$.
\end{proof}

\end{document}